\documentclass[11pt]{article}

\usepackage[all]{xy}

\font\tenmsb=msbm10
\font\sevenmsb=msbm7
\font\fivemsb=msbm5

\newfam\msbfam
\textfont\msbfam=\tenmsb
\scriptfont\msbfam=\sevenmsb
\scriptscriptfont\msbfam=\fivemsb
\def\Bbb#1{{\fam\msbfam #1}}

\font\teneufm=eufm10
\font\seveneufm=eufm7
\font\fiveeufm=eufm5
\newfam\eufmfam
\textfont\eufmfam=\teneufm
\scriptfont\eufmfam=\seveneufm
\scriptscriptfont\eufmfam=\fiveeufm

\newtheorem{tm}{Theorem}[subsection]
\newtheorem{lm}[tm]{Lemma}
\newtheorem{pr}[tm]{Proposition}
\newtheorem{rmk}[tm]{Remark}
\newtheorem{cor}[tm]{Corollary}

\newtheorem{??}[tm]{Question}
\newtheorem{defi}[tm]{Definition}

\newtheorem{ass}[tm]{Assumption}

\def\lorw{\longrightarrow}
\newcommand\n{\noindent}

\newcommand\ci{\cite}

\newcommand\rat{{\Bbb Q}}

\newcommand\zed{{\Bbb Z}}

\newcommand\pn[1]{{\Bbb P}^{#1}}

\newcommand\blacksquare{{\hspace*{\fill} $\fbox{}$}}

\newcommand{\fxn}[3]{{#1}_* \rat_{#2} [{#3}]  }

\newcommand{\im}{ \hbox{\rm Im \,} }

\newcommand{\ke}{ \hbox{\rm Ker \,} }

\newcommand{\td}[1]{ \tau_{ \leq {#1} } }

\newcommand{\tu}[1]{ \tau_{ \geq {#1} } }

\newcommand{\ih}[2]{  I\!H^{#1}({#2}) }

\title{Hodge-theoretic aspects of the Decomposition Theorem}

\author{Mark Andrea A. de Cataldo\thanks{Partially supported by N.S.F.
Grant DMS 0202321 and  0501020}$\;$  
Luca Migliorini\thanks{ Partially supported by
G.N.S.A.G.A.}
}

\date{February 5, 2006}

\begin{document}\maketitle

\begin{abstract}
Given a projective morphism of compact, complex, algebraic varieties
and a relatively ample line bundle on the domain
we prove that  a suitable choice, dictated by the line bundle,
of the decomposition isomorphism of the Decomposition Theorem
of Beilinson, Bernstein, Deligne and Gabber, yields
isomorphisms of pure Hodge structures.
The proof is based on a new cohomological characterization
of the decomposition isomorphism associated with the line bundle.
We prove some corollaries concerning the intersection form
in intersection cohomology, the natural map from cohomology 
to intersection cohomology, projectors and Hodge cycles,
and induced morphisms in intersection cohomology.
\end{abstract}

\tableofcontents

\section{Introduction}
\label{intro}
Let $f :X \to Y$ be a projective map of proper, complex, algebraic
varieties. The Decomposition Theorem predicts that the derived
direct image complex $Rf_* IC_X$ of the rational intersection cohomology 
complex $IC_X$ 
 of $X$ splits into the direct
sum of shifted intersection cohomology complexes on $Y.$
This splitting is {\em not} canonical. When viewed
in hypercohomology, it yields decompositions
of the rational intersection cohomology groups $I\!H(X,\rat)$ into
the direct sum of intersection cohomology groups with twisted coefficients
of closed subvarieties of $Y.$

The Decomposition Theorem is the deepest known fact concerning the homology
 of complex algebraic 
varieties and it has far-reaching consequences.
The following consideration may give a measure of the importance as well as 
of the special character of this result.
The splitting behavior of $Rf_* IC_X$ over $Y$ is dictated in part
by the one over any open subset $U \subseteq Y.$ 
This remarkable fact has no counterpart in other geometries, e.g.
complex analytic geometry, real algebraic geometry, etc.
More precisely:
let $U \subseteq Y$ be a Zariski-dense
open subset, $S
\subseteq U$ be a closed submanifold, ${\cal L}$ be a local system,
i.e. a locally constant sheaf, on $S;$
assume that a shift ${\cal L}[l]$ is a direct summand of 
$(Rf_* IC_X)_{|U}$ on $U;$ then a certain shift of 
the intersection cohomology complex
$IC_{\overline S}({\cal L})$
 on the closure $\overline{S} \subseteq Y$
 is a direct summand of $Rf_* IC_X$ on $Y.$

However,  the decomposition isomorphism is {\em not} canonical and
it is not clear, 
and in fact  not true, that
the various additional structures present in the various intersection 
cohomology  groups involved
should be preserved under an arbitrary splitting.
Let us consider the example of resolution of singularities.
In this case the Decomposition Theorem predicts the existence of
splitting injections $I\!H(Y, \rat) \to H(X,\rat).$
One may ask if it is possible to  realize $I\!H(Y,\rat) $
 as a sub-Hodge structure
of the pure Hodge structure $H(X, \rat).$
More generally, given any projective  map $f$, one may ask the 
same question for all the summands 
(see (\ref{byst}) and (\ref{dshtr})) of $I\!H(X,\rat)$
arising from the Decomposition 
Theorem. 

In this paper,  by building on our previous work \ci{htam},  we answer this question affirmatively in 
Theorem \ref{tuttoht} by checking that  a certain decomposition isomorphism 
 $g_{\eta},$ that depends on the choice of
an $f-$ample line bundle $\eta$ on $X,$  
 turns out to do the job.

Let us summarize the contents of this paper.
Given a projective map $f:X \to Y$ as above and an $f-$ample
line bundle $\eta$ on $X,$ Deligne, in \ci{shockwave} has constructed a distinguished
decomposition isomorphism
$g_{\eta}.$ 
Theorem \ref{tuttoht} shows that $g_{\eta}$ induces  an isomorphism
(\ref{dshtr}) of pure Hodge structures.
Let us emphasize that while  this isomorphism, being an isomorphism in the derived category,
is of a local nature, our result implies that it has global consequences, 
specifically, concerning the pure Hodge structure in intersection cohomology.
The proof is based on Proposition \ref{okt}, i.e. on
a property of $g_{\eta}$ expressed via its primitive components
$f_{i, \eta}.$

We obtain the following Hodge-theoretic consequences: Theorem \ref{c5},
on the intersection pairing on intersection cohomology;
Theorem \ref{alce}, on the natural map
$a_Y: H(Y,\rat) \to I\!H(Y,\rat);$ 
Theorem \ref{prcon},
on the homological cycles associated with the Decomposition 
Theorem;
Theorem \ref{imi}, on the morphism induced
by a surjective $f$ in intersection cohomology.

\bigskip
\n
{\bf Acknowledgments.} The first-named author
thanks the University of Bologna, the C.N.R.
of the Italian Government, the U.N.A.M
of Mexico City and the C.I.M.A.T. of Guanajuato for partial support.
The second-named author thanks  the N.S.F.
for partial support.

\subsection{Notation and preliminaries}
\label{notbm}
We fix the following notation. See also \ci{htam}.
For an introduction to the decomposition theorem with some examples worked out see
 \ci{leiden}.

-- $f :  X^n \to Y^m:$ a projective map of compact, complex, algebraic
varieties of the indicated dimensions.

-- $\eta: $  the first Chern class of a $f-$ample line bundle on $X.$

-- $f_{\bullet}= Rf_*:$ the derived direct image functor.

-- $H(X)= {\Bbb H} (X, \rat_X)= {\Bbb H}(Y, f_{\bullet} \rat_X):$
the $\rat-$singular cohomology of $X;$
at times we omit seemingly unnecessary cohomological degrees.

-- $IC_X:$ the intersection cohomology complex $X$ with $\rat-$
coefficients; if $X$ is smooth, then $IC_X=\rat_X[n].$

-- $I\!H^{n+l}(X)= {\Bbb H}^l(X, IC_X) = {\Bbb H}^l(Y, f_{\bullet}
IC_X), l \in \zed:$ the $\rat-$intersection cohomology groups of $X.$

-- $D_Y:$ the bounded derived category of constructible sheaves
on $Y$ of $\rat-$vector spaces, endowed with the $t-$structure
associated with middle-perversity.

-- $P_Y:$ the abelian category of perverse sheaves on $Y;$
$P_Y \subseteq D_Y$ is the heart of the middle-perversity $t-$structure.

-- $^p\!H^i: D_Y \to P_Y:$ the associated cohomological functors.


-- PHS, MHS, SHS: pure, mixed and   Hodge sub-structure.

If $a: K \to K'$ is a morphism in $D_Y,$ then we often use the same symbol
for the induced map in hypercohomology.

The category $P_Y$   is Artinian  and the Jordan-H\"older Theorem holds.
 The simple objects
are the intersection cohomology complexes $IC_{\overline{Z}}({\cal L})$
where $Z\subseteq Y$ is a smooth locally closed subvariety and ${\cal L}$
is a simple local system on $Z.$ A semisimple object of $P_Y$
is a finite direct sum of such objects.

The following results have been first proved by Beilinson,
Bernstein, Deligne and Gabber in  \ci{bbd}
using algebraic geometry in positive characteristic.
M. Saito has proved them
in the more general context of mixed Hodge modules in  \ci{samhp}.
We have given a proof in \ci{htam} using classical Hodge theory.
The earlier paper  
\ci{demigsemi}  had dealt with the special, but revealing case
 of semismall maps.

\begin{tm}
\label{dttm}
({\bf Decomposition Theorem (DT)})
There exists an isomorphism in the derived category $D_Y$
$$
\phi \; : \;   \bigoplus_{i \in \zed}{\,  ^p\!H^i( f_{\bullet} IC_X) [-i] }  
\; \simeq \; f_{\bullet}IC_X , \qquad
^p\! H^i( f_{\bullet} IC_X) \; \mbox{ semisimple in } \; P_Y.
$$
\end{tm}

The Chern class $\eta$  defines
a map 
$
\eta: IC_X \to IC_X[2].
$ 
Taking push-forwards and cohomology we get
maps
$$
e \;  :=\;  ^p\!H^j(f_{\bullet} \eta) \; : \;  ^p\!H^j (f_{\bullet} IC_X) 
\lorw \,  ^p\!H^{j+2}(f_{\bullet}IC_X).
$$

\begin{tm}
\label{rhltm}
({\bf  Relative Hard Lefschetz  Theorem}) 
For all $i\geq 0$ the map
$$
e^i \; :\; 
^p\! H^{-i}(f_{\bullet} IC_X ) \; 
\to \;
 ^p\!H^i(f_{\bullet} IC_X), 
$$
is an isomorphism.
\end{tm}

Let us collect together 
some well-known facts that we shall use.

Let $Y$ be a proper algebraic variety. Goresky-MacPherson defined 
the intersection homology
using a sub-complex of the complex of geometric chains of $Y.$ This gives rise
to a natural map $I\!H_l(Y) \to H_l(Y).$
Using the perfect pairing in intersection (co)homology mentioned below,
 there is the  natural dual  map
$H^l(Y) \to I\!H^l(Y)=I\!H_{2m-l}(Y)=I\!H_l(Y)^{\vee}.$  
This map can be described also
as the map in hypercohomology stemming from the natural map
$a_Y: \rat_Y[m] \to IC_Y$ that corresponds to $1 \in \rat$ under the isomorphism
\begin{equation}
\label{qic}
\rat = H^0(Y)=  {\rm Hom}(\rat_Y[m], \rat_Y[m]) \simeq 
{\rm Hom}( \rat_Y[m], IC_Y) =
I\!H^0(Y) = \rat
\end{equation}

The number $1 \in \rat= H^0(X),$ ${\rm Id}$ and the adjunction map $
adj(f) : \rat_Y \to f_{\bullet} \rat_X $
correspond to each other under the isomorphisms
\begin{equation}
\label{cit}
\rat = H^0(X) = {\rm Hom}(\rat_X , \rat_X ) \simeq 
{\rm Hom} (\rat_Y,  f_{\bullet} \rat_X ).
\end{equation}
 
The map $adj(f)$ induces the familiar pull-back
in cohomology $f^*: H(Y) \to H(X).$

By adjunction and by (\ref{qic}) applied to $X:$
 \begin{equation}
\label{rana} 
{\rm Hom}( \rat_Y[n], f_{\bullet}IC_X)= {\rm Hom}(\rat_X[n], IC_X)  \simeq \rat. 
\end{equation}

\begin{rmk}
\label{ranami}
{\rm
The equalities above  hold if we replace
$Y$ by a {\em connected} open subset $U \subseteq Y$
and $X$ by $f^{-1}(U).$
}
\end{rmk}

Given a  proper variety $Y,$ there is a non-degenerate
 intersection pairing on intersection cohomology:
$$
I\!H^{n-l}(Y) \times I\!H^{n+l}(Y) \lorw \rat. 
$$
It has been first defined by Goresky-MacPherson in terms of geometric cycles.
It affords an alternative description as the map in hypercohomology
stemming from the duality isomorphism 
$$
d_Y \; :  \; IC_Y \; \simeq  \; IC_Y^{\vee},
$$
On the regular part, the isomorphism $d_Y$ 
coincides with the usual Poincar\'e
duality isomorphism. One way to say this is the following:
${IC}_{Y_{reg}}$ is canonically isomorphic to $\rat_{Y_{reg}}[m]$
and the duality isomorphism for $IC_Y$ is the unique morphism
in ${\rm Hom}(IC_{Y}, IC_Y^{\vee}) \simeq \rat$ extending
the duality isomorphism for $\rat_{Y_{reg}}[m].$

\subsection{Review of \ci{htam}}
\label{totrec}
We recall some of the result of our paper
\ci{htam} in the form we need them here.

For every $l \in \zed,$ $I\!H^l(X)$ carries a canonical
PHS of weight $l.$

The subspaces of the perverse filtration
(\ref{filtr}) $I\!H^l_{\leq i}(X) \subseteq
I\!H^l(X),$ $i \in \zed,$ are SHS.
In fact,  the filtration $I\!H^l_{\leq i}(X)$ can be described 
up to shift as the  monodromy weight filtration of the endomorphism $I\!H(X) \to I\!H(X)$
given by the cup-product with the pull back of {\em any} ample bundle on $Y.$
The graded pieces (\ref{pco}) $I\!H^l_i(X)= 
I\!H^l_{\leq i}(X)/I\!H^l_{\leq i-1}(X)$
inherit the quotient PHS.

Let $Y= \coprod_{d\geq 0}{ S_d} = \coprod_{d \geq 0} \coprod{S},$
be a stratification of $Y$ adapted to $f,$  where $S$ ranges over
the connected components of the $d-$dimensional stratum $S_d.$
There is a canonical decomposition  given  by strata for the semisimple:
 \begin{equation}
 \label{sdec}
 ^p\!H^i( f_{\bullet} IC_X) 
\;
= \; \bigoplus_{d \geq 0   }\bigoplus_{S \subseteq S_d}{IC_{\overline{S}}
 ( {\cal L}_{i,S}   )},
 \end{equation}
 where ${\cal L}_{i,S}$ are semisimple local systems on $S.$
 
 The ensuing decomposition in hypercohomology is by SHS:
 \begin{equation}
\label{byst}
I\!H_i^{n+l}(X) \; = \; 
{\Bbb H}^l(Y,  \, ^p\!H^i( f_{\bullet} IC_X) [-i]    ) \; = \; 
 \bigoplus_{d,S}{ {\Bbb H}^l(Y, IC_{\overline{S}} ( {\cal L}_{i,S}   )[-i] ) 
},
 \qquad \forall i,l \in \zed,
\end{equation}
where the first equality stems from (\ref{pco}).

There are Hard Lefschetz isomorphisms for the action of  $\eta$ 
on the graded pieces:
\begin{equation}
\label{hletal}
e^i \; : \; I\!H^l_{-i}(X) \simeq I\!H^{l+2i}_{i}(X), \quad \forall l \in \zed, \; \forall i\geq 0.
\end{equation}

A natural question, see \ci{bbd} and \ci{mac83}, is whether
the decomposition
\begin{equation}
\label{dshtr}
\bigoplus_i \,  \phi \,   ( I\!H^{n+l}_i (X)       )\; =\;  I\!H^{n+l}(X),
\end{equation}
its refinements stemming from (\ref{byst})
and the further refinements stemming from the $(\eta,L)-$de\-composition
we prove in \ci{htam},
 are  isomorphisms of PHS
for a suitable choice of the isomorphism 
$\phi.$

Our main Theorem \ref{tuttoht} gives  a positive answer.

We shall need the following simple

\begin{lm}
\label{bho}
Let $A$ and $B$ be  rational Hodge structures and
$$
A \stackrel{g}\lorw B \stackrel{p}\lorw A
$$
 be  linear maps with $p\circ g =Id,$ and $p$  be a map  of rational Hodge structures
 and $g(A) \subseteq B$ a SHS.
 
 \n
 Then $g$ is a map of HS.
\end{lm}
{\em Proof.} We need to show that, after complexification, $g(A^{pq}) \subseteq B^{pq}.$
Let $a_{pq} \in A^{pq}.$ 
We have that $g(a_{pq})= \sum{ b_{st} }$ for unique $b_{st} \in B^{st}.$ 
Noting that $g$ is necessarily injective and since we are assuming that 
$g(A) = \oplus{ (g(A) \cap B^{pq})},$
then $b_{st} =g(c_{st})$ for a unique collection  $c_{st} \in A.$
Since $a_{pq} = \sum{ p(b_{st})},$ we have that
$p( b_{st}) =0$ for $(s,t) \neq (p,q)$ and we also  have that $0= p(g(c_{st})) = c_{st}$
for the same $(s,t).$ It follows that $a_{pq}= c_{pq}$ and that
$g(a_{pq}) = g(c_{pq}) = b_{pq}.$
\blacksquare

\begin{rmk}
\label{silly}
{\rm
The example of $A =B$ as vector spaces, but with conjugate Hodge structures,
shows that having $g$ defined over $\rat$ and having image a SHS is not
sufficient to have a map of HS.
}
\end{rmk}

\section{Formalism for decompositions}
\label{prelim}
The aim of this paper is to show that the  isomorphism $g_{\eta}, $ constructed by Deligne
in \ci{shockwave},
gives rise to an isomorphism of PHS.

The morphism $g_\eta$ is constructed by assembling together
certain morphisms $f_{i, \eta}$ defined on the primitive components
$P^{-i}_{\eta} \subseteq 
 \, ^p\!H^{*}(f_{\bullet} IC_X)$
for the action of $\eta$ on the perverse cohomology
complexes $^p\!H (f_{\bullet} IC_X).$

In this section we review the constructions of $g_{\eta}$ and $f_{i,\eta}$
given in \ci{shockwave}.
We  then prove Proposition \ref{okt} that is the key
to our main result Theorem \ref{tuttoht}.

To simplify the notation, we present most of the material of this section
in the abstract context of a triangulated category
with $t-$structure.
For our purposes,
the main example of the formalism discussed below is
given by ${\cal D}= D_Y,$ $ K= f_{\bullet} IC_X,$ $
F(-) = {\Bbb H}^0(Y,-),$ etc.

A geometric example is discussed in $\S$\ref{exqu}.

\subsection{Decomposition via $E_2-$degeneration}
Let ${\cal D}$ be a triangulated category with $t-$structure. 
Its heart ${\cal C} \subseteq {\cal D}  $ 
is an abelian category. This data comes equipped with
the corresponding
cohomology functors  $H^i: {\cal D} \to {\cal C}.$

We consider  objects  $K$ of ${\cal D}$
with bounded cohomological amplitude,  i.e. such that
$H^i(K)  =0, $   for every $ |i| \gg 0.$
We also assume the $t-$structure non-degenerate, see, \ci{bbd}, 1.7.
This implies that  $H^i(K)=0$ for all $i$
if and only if $K=0.$

For any  object $X$
of ${\cal D}$  there is a spectral sequence
\begin{equation}
\label{ss1}
E_2^{pq}={\rm Hom}(X[-p],H^q(K)) \Longrightarrow  {\rm Hom}(X[-p],K[q]).
\end{equation}

If we assume that (\ref{ss1}) is  $E_2-$degenerate for $X=H^i(K)$, 
for any $i$,  then there exists
an isomorphism in ${\cal D}:$
\begin{equation}
\label{ptp}
\phi\; := \; \sum_i  \phi_i\; : \; 
\bigoplus_i H^i(K)[-i] \; \stackrel{\simeq}\lorw \; K
\end{equation}
inducing the identity in cohomology.
This can be seen as follows.
The $E_2-$degeneration ensures  that
${\rm Hom}(H^i(K), H^i(K))$ is a quotient of ${\rm Hom}(H^i(K)[-i], K).$  
This implies that
every map $H^i(K) \to H^i(K)$ admits
a, {\em not necessarily unique}, 
 lift  to ${\rm Hom}(H^i(K)[-i],K).$ By applying this to ${\rm Id}: H^i(K) \simeq H^i(K,)$ we get a map
$\phi_i: H^i(K)[-i] \to K$ inducing the identity in cohomological degree $i$
and the zero map in the remaining degrees.
By putting together these maps, we get  the morphism (\ref{ptp}) that, being
 the identity  in cohomology,   is an isomorphism in $D.$

Any isomorphism $\phi$ as above can be normalized by
an automorphism of the lhs so that it induces the identity
in cohomology. We always work with such
 isomorphisms.

 \medskip
In short, the degeneration
of (\ref{ss1}) implies the  existence of a splitting (\ref{ptp}).
 However, as the construction shows,
 this decomposition
is { \em not unique.}

\subsection{$E_2-$degeneration via the Deligne-Lefschetz Criterion}
\label{dlcrt}
Let $F: {\cal D} \to Ab$ be a cohomological functor. As usual, 
set $F^0(X):= F(X)$
and $F^l(X): = F^0(X[l]).$
Fix a morphism
\begin{equation}
\label{eta}
 \eta: K \to K [2].
\end{equation}  
For $a\in F^l(K),$ denote $\eta (a)$ by $\eta \wedge a \in F^{l+2}(K).$ 
Set  $e: = H^l (\eta):H^l(K) \to H^{l+2}(K).$

\begin{ass}
\label{asshle}
Assume that $\eta$ satisfies the following Hard Lefschetz relation:
\begin{equation}
\label{hle}
 e^i : H^{-i}(K) \simeq H^i(K), \qquad \forall i\geq 0.
\end{equation}
\end{ass}

The Deligne Lefschetz Criterion (cf. \ci{dess}
and \ci{shockwave}, p.116) is a sufficient condition 
for degeneration and splitting:
 the Hard Lefschetz relation  (\ref{hle}) implies that
 the spectral sequence (\ref{ss1}) is $E_2-$degenerate so that
there exist splittings  $\phi$ as  in (\ref{ptp}).

\bigskip
The main example for us is the following. 
Let $f: X \to Y$ be  a  projective morphism  of varieties,   $\eta
\in  {\rm Hom} (\rat_X, \rat_X [2])$  be the first Chern class
of an $f-$ample  line bundle on $X.$ 
Setting   $K: = f_* IC_X.$ we have $f_*\eta: K \to K[2]$ etc.
 The Relative Hard Lefschetz Theorem
 \ref{rhltm} holds and one deduces from it
the Decomposition Theorem \ref{dttm} (without the semisimplicity
assertion). 

 \subsection{Primitive Decomposition}
 \label{prde}
Since  the heart ${\cal C}$ of  the given $t-$structure 
on ${\cal D}$ 
is an abelian category,   with slight abuse of 
 language, we think of kernels and images in ${\cal C}$  as subobjects. 
 
By analogy with the classical  primitive decomposition of the cohomology
of a projective manifold with respect
to an  ample line bundle we define: 
$$
P^{-i}_{\eta}: =  \ke{ \{ e^{i+1}: H^{-i}(K) 
\to H^{i+2}(K)} \}, \qquad  i\geq 0,
$$

$$
e^jP_{\eta}^{-i} :=  \im{  \{ e^j: P_{\eta}^{-i} \to H^{2j-i}(K)\}  }, \qquad 0 \leq j \leq i.
$$

There is the   Lefschetz-type
canonical decomposition isomorphism in the heart ${\cal C}:$

\begin{equation}
\label{ldo}
\bigoplus_{l = 2j-i; \;
0 \leq j \leq i}{ e^jP_{\eta}^{-i} }\;  \simeq \;  
H^{l}(K).
\end{equation}


\subsection{The $t-$filtration}
\label{tfil}
Let $F: {\cal D} \to Ab$ be a  cohomological functor. The $t-$structure
on ${\cal D}$ defines a filtration on  the groups $F^l(K):$
\begin{equation}
\label{filtr}
F^l_{\leq i} (K) \;  := \;  \im{  \{  F^l (\td{i} K)  \to F^l (K ) \}  }.
\end{equation} 
This filtration is the abutment of the spectral sequence (\ref{ss1})
and we call it 
the $t-$filtration. 

In the geometric case,
we get an increasing filtration 
$I\!H_{\leq i} (X) \subseteq I\!H(X)$ and  we call it the {\em
 perverse filtration}.

For every isomorphism
$
\phi\; := \; \sum_i  \phi_i\; : \; 
\bigoplus_i H^i(K)[-i] \; \simeq \; K
$
we have
$$
 F^l_{\leq i} (K)  \; = \; F^l(\phi( \bigoplus_{i' \leq i}     H^{ i'}(K)[-i'])  ;
$$
 this means that, while the individual summands
on the rhs are not canonically embeddable in the lhs,
the images of the  direct sums above are the canonical subspaces  yielding
the $t-$filtration.

By abuse of notation we often denote with the same symbol a map of, say, 
complexes and the resulting map in, say, hypercohomology. 
 
Since $K$ decomposes, the associated graded pieces satisfy canonically
 \begin{equation}
 \label{pco}
 F^l_i(K) \;  : \;  = F^l_{\leq i} (K) / F^l_{\leq i-1}(K) 
\; \simeq \;  F^l ( H^i(K)[-i]), \qquad
 \forall i, l \in \zed.
 \end{equation}

Since $\eta$ (\ref{eta})  is a $2-$morphism, 
we have
$$
\eta^j \; :   \;  F^l_{\leq i}(K) \lorw F^{l+2j}_{\leq i+2j} (K), 
\qquad \forall i, l \in \zed, \;
\forall j \geq 0.
$$
For every $i \geq 0,$ the composition 
$$ 
 F^{l } (   H^{-i}(K)[i]  ) \;   \stackrel{\phi_|}\to (  F^{l } (K) ) )  
  \stackrel{\eta^i}\lorw {F}^{l}(K[2i]) \stackrel{pr_i \circ \phi^{-1}}\lorw 
F^{l} (   H^i (K)[i]  )
$$
coincides, in view of (\ref{pco}) and (\ref{hle}), with the isomorphism
(see (\ref{hletal})):
$$ 
e^i \; : \;  F^l_{-i}(K ) \;  \simeq \;  F^{l+2i}_i (K).
$$
and the composition
\begin{equation}
\label{isonim}
\phi (  {F}^{l} (  P_{\eta}^{-i-2j} (-j)  [i]  ) )  
\stackrel{\eta^{i+1}}\lorw {F}^{l +2(i+1)}_{\leq i+2}(K) 
\lorw 
{F}^{l+2(i+1)}_{i+2}(K)
\end{equation}
is an isomorphism onto its image for every  $j>0.$

  The condition  $a \in F^l( P_{\eta}^{-i}[i]),$  does not imply that 
 $\phi(a) $ is  primitive in the usual sense, i.e. $\eta^{i+1} \wedge \phi(a)=0.$
What is true is (\ref{cvdep}) below.
  An element $a \in F^l( P_{\eta}^{-i}[i]),$ 
$i \geq 0,$  satisfies $e^{i+1}a =0 \in F^l (H^{i+2}(K)[i]  ) =
  F^{l+2(i+1)} ( H^{i+2}K [-i-2]).$ This means that for every splitting  $\phi$
  as in  (\ref{ptp}) we have that $\phi(a) \in F^l_{\leq -i}(K)$ and 
  \begin{equation}
  \label{cvdep}
\eta^{i+1}  \wedge \phi (a)  \in F^{l+2(i+1)}_{\leq  i+1   }(K) 
\subseteq F^{l+2(i+1)}_{\leq  i+2   }(K)
\end{equation}
so that $\eta^{i+1} \wedge \phi(a)$ becomes zero in $F^{l+2(i+1)}_{i+2}( K).$

 \subsection{The canonical morphisms $f_{i,\eta}: P_{\eta}^{-i}[i] \lorw K$}
 \label{cafi}
For the reader's convenience,
we now recall Deligne's  construction of the maps $f_{i, \eta}$ (cfr.\ci{shockwave}).

We assume \ref{asshle}. Since  then $K$ splits, for every cohomological functor $F: {\cal D}
\to Ab$ we have short exact sequences
$$
0 \lorw F^{i+1} H^{-i-1}K \lorw F^0 \tu{-i-1} K \lorw
F^0 \tu{-i} K \lorw 0.
$$

Let $i \geq 0$ and $t: {\cal T} \to H^{-i}K$
be a morphism in ${\cal D}$ that factors 
through $P^{-i}_{\eta}.$
In particular, we have $0=\eta^s \circ t:
{\cal T} \to H^{-i +2s } K,$ for every $s >i.$
The morphism $t$ induces a morphism $x: {\cal T}[i] \to \tu{-i}K.$
Let $T:= {\rm Hom} ({\cal T}[i] ,-) : {\cal D} \to Ab.$

\begin{pr}
\label{so}
Let $t \in T^0 H^{-i}K[i]$ and $x \in T^0 \tu{-i}K$ be as above.
There exists a unique lift $F^{-i}K \ni \tau:
{\cal T}[i] \to K$
of $x$ such that $0= \eta^s \circ \tau\in
T^{2s} \tu{s} K,$ for every $s >i.$
\end{pr}
{\em Proof.} (See \ci{shockwave}, Lemme 2.2).
There is the commutative diagram of 
 short exact sequences:
 $$
 \xymatrix{
 0  \ar[r]      & T^{i+1}  H^{-i-1} K   \ar[r] 
  \ar[d]^{e^{i+1}}_{\simeq} & T^0   \tu{-i-1} K    
 \ar[r]   \ar[d]^{ \eta^{i+1} }  & T^0  \tu{-i} K    
 \ar[r]  \ar[d]^{\eta^{i+1}}  &  0 
 \\
   0 \ar[r]  & T^{i+1}  H^{+i+1} K(i+1)  \ar[r]  & T^{2(i+1)}  
\tu{i+1} K    
 \ar[r]  & T^{2(i+1)}  \tu{i+2} K    \ar[r]   &  0  
  }
 $$
Since $\eta^{i+1}\circ x = 0$
and $\eta^{i+1}$ is an isomorphism, the Snake Lemma yields
 the existence of a unique lift of $x,$ 
$x_{-i-1} \in T^0 \tu{-i-1} K,$
with the property
that
$\eta^{i+1} \circ x_{-i-1}=0.$
Repeating this procedure, with $i$ replaced by $i+1$ and
$x$ by $x_{-i-1},$ that clearly satisfies $\eta^{i+2} \circ x_{-1-2}=0,$
we obtain, for some $r \gg 0,$  $\tau:= x_{-r} \in T^0\tu{-r} K = T^0 K = 
{\rm Hom} ( {\cal T}[i], K    )$ with the required property.
\blacksquare

\begin{defi}
\label{deffi}
{\rm
Fix $i \geq 0.$ Let $t: T:=P^{-i}_{\eta} \to H^{-i}(K)$
be the inclusion. Proposition \ref{so} yields morphisms:
$$ 
f_{i, \eta} \; : \; P^{-i}_{\eta} \lorw K.
$$
}
\end{defi}

These morphisms are
characterized by the two properties  that

\smallskip
$(i)$ $H^{-i}(f_{i, \eta}):  P_{\eta}^{-i} \to H^{-i}(K)$ 
is the natural inclusion and

$(ii)$ for  every $s >i,$ the composition below is zero:
$$
P_{\eta}^{-i} [i] \lorw K \stackrel{\eta^s}\lorw  
K[2s]  \lorw (\tu{s} K)  [2s].
$$

 The second condition implies
that if $F: {\cal D} \to Ab$ is any cohomological
functor and $\phi$ is any decomposition isomorphism
 (\ref{ptp}) coinciding with $f_{i,\eta}$ 
on  the summand $P^{-i}_{\eta} [i],$  then we have 
\begin{equation}
\label{sw}
\eta^{s}  \wedge  \phi (   F^l( P^{-i}[i]  ) ) \; 
\subseteq \; F^{l+2s}_{\leq s -1} (K) ),
\qquad \forall s>i.
\end{equation}

By (\ref{cvdep}), 
a priori 
the lhs is contained in  the bigger $F^{l+2s}_{\leq s-1 + (s-i)} (K).$
This is an important restriction and is used in our proof of
the key Proposition \ref{okt}.
We shall discuss it further in a geometric example in $\S$\ref{exqu}.

\begin{rmk}
\label{fdepe}
{\rm
The objects $P_{\eta}^{-i},$ depend on $\eta$ and so do the
morphisms $f_{i, \eta}.$ It may happen
that
$P^{-i}_{\eta}$ is independent of $\eta.$
It is important to keep in mind that, even in this case,
the morphisms $f_{i, \eta}$ depend on $\eta.$ See the example
of $\S$\ref{exqu}.
This explains why in general one cannot hope for a canonical decomposition
isomorphism (\ref{ptp}). Of course, in special cases, one may have a
 distinguished choice of $\eta$
and consider the resulting
$g_{\eta}$ canonical.
}
\end{rmk}

\subsection{The isomorphism $g_{\eta}:  \oplus H^{l}(K) [-l] \simeq K.$}
\label{deigo}
We assume \ref{asshle} and  therefore we have the 
morphisms $f_{i,\eta}$ of $\S$\ref{cafi}.

The isomorphism $g_{\eta}$  is constructed
by assembling together the $f_{i, \eta}$ using the primitive
Lefschetz  decomposition (\ref{ldo}).

We start by defining $g_{l, \eta}: H^l K [-l] \to K$
by first defining it
on the direct summands $e^jP_{\eta}^{-i} [i-2j],$
$ 0 \leq j \leq i, $ $l = 2j-i,$ 
as the composition
$$
g_{l, \eta} \; : \; e^jP_{\eta}^{-i} [i -2j] \stackrel{ (e^j)^{-1}}\lorw
P^{-i}_{\eta} [i -2j] \stackrel{f_{i,\eta}[i-2j]     }\lorw
K[-2j] \stackrel{\eta^j}\lorw K. 
$$ 
Collecting together the maps $g_{l, \eta},$ $l \in \zed,$ we obtain
a decomposition  isomorphism 
\begin{equation}
\label{deliso}
g_{\eta} \; : \; \bigoplus_l H^l(K)[-l] \; \simeq \; K.
\end{equation}
It depends on $\eta:$
the $g_{l, \eta}$ are obtained
via the $f_{i,\eta}$ and through repeated applications of $\eta.$
It induces the identity in cohomology and,
by construction, the restriction of $g_{\eta}$ to the direct summand 
$P_{\eta}^{-i}[i]$
is $f_{i,\eta}.$

The properties of $g$ which are relevant to this paper are
the following.

Let $0 \leq j \leq i.$ For every $j'$ s.t. $j'+j \leq i,$ 
we have that $g^{-1}_{\eta} \circ \eta^{j'} \circ g_{\eta} =
 e^{j'}$ when restricted to  $e^jP^{-i} [i-2j].$
 In particular, the cup product with $\eta^{j'}$
has the simplest possible expression in terms of the direct sum decomposition,  i.e.
the following diagram is commutative:
$$
\xymatrix{
F^{l}   ( e^jP^{-i} [i-2j]  )    \ar[r]^{g_{\eta}}   \ar[d]^{e^{j'} } & 
g (  {F}^{l }   ( K ) )           \ar[d]^{\wedge {\eta}^{j'} }
\\
F^{l +2j'}  ( e^{-j'-j}P^{-i}[i-2j] )    \ar[r]^{g_{\eta}}    & 
g  ( F^{l +2j'} ( K )  ) ,
}
$$
or, in words,  $\eta^{j'}$ and $g$ commute when applied to
elements of the primitive decomposition as long as $\eta^{j'}$ stays
injective as predicted by the Hard Lefschetz  property \ref{asshle}.

In the remaining range, we have the key restriction
 (\ref{sw}).

There is the   decomposition  
\begin{equation}
\label{delde}
F^l(K) =  \bigoplus_{ 0 \leq j \leq i}{ \eta^j \wedge  
f_{i, \eta} ( F^{l-2j} ( e^jP_{\eta}^{-i}[i]  )) },
\end{equation}
i.e. the lhs can be built inductively using the images of primitives
via the maps $f_{i,\eta}$ 
in degrees $\leq l$ via cup products with $\eta.$

In our geometric situation, $K= f_{\bullet} IC_X,$ $F= {\Bbb H}^0(Y,-)$ etc, 
we get
$$
I\!H^{n+l} (X) = \bigoplus_{0 \leq j \leq i}{
\eta^j \wedge f_{i, \eta} 
(
{\Bbb H}^{l-2j} (Y, P^{-i}_{\eta} [i]  )
), \qquad l \in \zed.
}
$$

\subsection{ Characterization of 
$ f_{i,\eta}  (  {\Bbb H}^l (Y, P^{-i}_\eta [i] ) \subseteq
I\!H^{n+l}(X)$
}
\label{charoff}
We revert to our geometric situation: $K:= f_{\bullet} IC_X,$
$F(-):= {\Bbb H}^0(Y, -),$
etc. 

Fix $i \geq 0.$  We shall define maps $\Psi_t$ and express the images in 
hypercohomology 
$$
f_{i, \eta} ( {\Bbb H}^l ( Y, P^{-i}_{\eta} [i]    )   )
\subseteq 
g_{\eta} (   {\Bbb H}^l ( Y, \, ^p\!H^{-i}( f_{\bullet} IC_X   ) [i] ) ) 
\subseteq I\!H^{n+l}_{\leq -i} (X) \subseteq I\!H^{n+l}(X)
$$
as $\ke{ \Psi_{r-i}   },$ where  $r= r(f_{\bullet} IC_X)$ is
 the cohomological amplitude of  $f_{\bullet} IC_X.$
This will be  achieved by means of a repeated application
of the key restriction (\ref{sw}).

Let $g_{\eta}$ be the isomorphism  (\ref{deliso}) 
associated with $\eta.$

 In what follows,
for simplicity,
we omit some
cohomological degrees .

Consider the composition
$$
\Psi_0 \; : \; I\!H_{\leq -i}^{\bullet} (X) \stackrel{\eta^{i+1}}\lorw
I\!H^{\bullet+2(i+1)}_{\leq i+2} (X)
\lorw 
 I\!H^{\bullet+2(i+1)}_{ i+2} (X),
$$
and define inductively, for $t \geq 1:$
$$
\Psi_t \; : \; 
\ke{ \Psi_{t-1} } \stackrel{\eta^{i+t}}\lorw
I\!H^{\bullet+2(i+t)}_{ i+t} (X).
$$

\begin{pr}
\label{okt}
$$
\ke{\Psi_{r-i}}=  f_{i,\eta} ( {\Bbb H}^l(Y,   (P^{-i}_{\eta} [i] )  ) ).
$$
\end{pr}
{\em Proof.} We are going to prove by induction on $t$  that 
\begin{equation}
\label{isps}
\ke{\Psi_t}=  I\!H^{n+l}_{ \leq -i-t-1}(X)  \,    \oplus  \, 
f_{i, \eta} ( {\Bbb H}^l(Y,(P_{\eta}^{-i} [i] ) ) ), \qquad \forall t\geq 0.
\end{equation}
Taking $t  = r-i,$ where $r$ is the cohomological amplitude of $f_{\bullet}
IC_X,$  we can draw 
the desired conclusion, for $I\!H_{\leq -r-1}(X) =0.$

\n
We first prove (\ref{isps}) for $t=0.$
 We have
$$
{I\!H}^{n+l}_{\leq -i}(X)  =  I\!H^{n+l}_{\leq -i-1}(X)   \oplus
 f_{i, \eta} ( {\Bbb H}^l(Y,(P_{\eta}^{-i} [i] )) ) \oplus
\bigoplus_{j>0}{
 f_{i, \eta} ( {\Bbb H}^l(Y,
 (   P_\eta}^{-i -2j} (-j) [i]     )    ))     .
$$
The first summand lands first in  $ I\!H^{\bullet+2(i+1)}_{\leq i+1}(X)$
and is therefore  in the kernel of $\Psi_0.$
So is the second summand  since it  lands first in 
$I\!H^{\bullet+2(i+1)}_{\leq i}(X)$ by virtue of (\ref{sw}).
As to the third summand, it maps isomorphically to its image via $\Psi_0$ by
 (\ref{isonim}).
This proves the case $t=0.$

\n
Assume we have proved (\ref{isps}) for $t-1$ and let us prove it for $t.$
We have the composition
$$
\Psi_t: \ke{\Psi_{t-1}} \stackrel{\eta^{i+t}}{\lorw}
I\!H^{\bullet+ 2 (i+t)} _{\leq i+t}(X) \lorw I\!H^{\bullet+2(i+t)}_{i+t}(X)
$$
where, by the inductive hypothesis:
$$
\ke{\Psi_{t-1}}=  I\!H_{\leq -i-t-1 }^{n+l} (X)   \, \oplus \,
 g_{\eta}
(   {\Bbb H}^l (Y,  \,  ^p\!H^{-i-t}(f_{\bullet} IC_X)[i+t]    ) ) \, 
  \oplus  \,  f_{i, \eta} ( {\Bbb H}^l(Y,(P_{\eta}^{-i} [i] )) ).
 $$
We conclude as in the case $t=0.$
\blacksquare

\medskip
Using orthogonality with respect to the intersection pairing in $I\!H(X),$
we can re-word Proposition \ref{okt} as

\begin{cor}
\label{ocond}
$$
 f_{i, \eta} ( {\Bbb H}^l(Y,(P_{\eta}^{-i} [i] )) )
= I\!H^{n+l}_{\leq -i} \cap \cap_{t \geq 1}
(  \eta^{i+t} I\!H^{n+l -2(i+t)}_{-i-t}(X)    )^{\perp}.
$$
\end{cor}
This formula shows that the realization of intersection
cohomology as a sub-Hodge structure of the cohomology of a 
resolutions of isolated singularities of threefolds and fourfolds worked out in \ci{leiden} coincides with the one defined by $g_{\eta}.$ 

\subsection{The isomorphism $g_{\eta}$ is Hodge-theoretic}
\label{htdt}
For simplicity let us now assume that $X$ is smooth and
let us briefly discuss the PHS on the graded spaces $H^l_i(X).$
In the paper \ci{htam} we have identified, up to some trivial shifting 
procedure, the perverse filtration
$H^l_{\leq i}(X)$ arising from a map $f: X \to Y$
with the filtration associated  with the nilpotent action
on $H^*(X)$ of the first Chern class of an ample line bundle on $Y.$
Since this action is of type $(1,1),$ the filtration
is  given by SHS. Accordingly, the subspaces of the filtration
are  PHS and the graded pieces, $H^l_i(X)$ 
inherit canonical PHS.

The Decomposition Theorem does not ensure that the resulting decomposition
 $H^l(X) = \oplus \phi (H^l_i(X)  ) $
(\ref{dshtr}) into the sum of the graded pieces
can be made into an isomorphism of PHS.

We are about to prove that this is achieved by the isomorphisms
$g_{\eta}.$

\begin{tm}
\label{tuttoht} 
Let $f: X \to Y$ be a projective morphism of compact varieties,
$\eta$ be an $f-$ample line bundle on $X.$ Then  
$g_{\eta}$ induces isomorphisms of weight $l$ PHS
$$
g_{\eta} \; : \; 
 \bigoplus_{i}{ I\! H^l_i (X)  } \; \simeq \;
 I\!H^l(X) .
$$
\end{tm}
{\em Proof.} By Lemma \ref{bho}, it is enough to show that
$g_{\eta} ( I\!H^l_i(X)   ) \subseteq I\!H^l(X)$ is a SHS
for every $i \in \zed.$\

\n
The cup product map $\eta: I\!H(X) \to I\!H^{+2}(X)$
is a map of PHS.

\n
By virtue of the $\eta-$decomposition
(\ref{delde}) associated with  $g_{\eta},$ it is enough 
to show that each subspace 
$f_{i, \eta} ( {\Bbb H}^{l-n}
(Y,  P^{-1}_{\eta} [i] ) ))$ 
is a SHS.
This follows from Proposition \ref{okt} that exhibits
those subspaces as iterated kernels of maps of PHS.
\blacksquare

\begin{cor}
\label{ichtr}
Assume, in addition, that $f: X \to Y$ is a resolution of singularities. Then 
$IC_Y \subseteq P^0_{\eta} \subseteq  \, ^p\!H^0 ( f_{\bullet} IC_X)$
and the induced map
$$
f_{0,\eta} \; : \; I\!H(Y) \lorw H(X)
$$
is an injection of PHS.
\end{cor}
{\em Proof.} The inclusion 
$IC_Y \subseteq  \,  ^p\!H^0 ( f_{\bullet} IC_X)$
holds over the smooth part of $Y$ and the Decomposition Theorem
implies that the inclusion must hold over $Y.$
Since the complexes $^p\!H^{\neq 0} ( f_{\bullet} IC_X )$
are supported on a proper subvariety of $Y,$ the simplicity of
$IC_Y$ implies the inclusion $IC_Y \subseteq P^0_{\eta}.$

\n
The summand $I\!H(Y) \subseteq I\!H_0 (X)$ corresponds
to the dense stratum in the strata-like decomposition
(\ref{byst}) and is therefore a SHS.
We conclude by Theorem \ref{tuttoht}.
\blacksquare

\subsection{An example: the blow up of a quadric cone}
\label{exqu}
Let $f: X \to Y$ be the blowing up at the vertex $v \in  Y$
of the projective cone $Y$ over a quadric surface
$\pn{1} \times \pn{1} \simeq Q \subseteq \pn{3}.$
There is the canonical $\pn{1}-$bundle projection
$p:X \to Q$ with sections $D:=f^{-1}(v)$ and 
$D_{\infty} := f^{-1}(\Delta_{\infty} ),$
where $\Delta_{\infty}\subseteq Y$ is the quadric at infinity.
There are the two surfaces $D_i: = p^{-1} (l_i),$
$i=1,2,$ where $l_i$ are two lines of the two distinct rulings of $Q.$
Let ${\Delta_i}:= f(D_i).$

We have the following relations in the $3-$dimensional group $H^2(X):$
$$
H^2(X) = \langle  D, D_1, D_2  \rangle 
= \langle D_{\infty}, D_1, D_2 \rangle, 
\qquad
D_{\infty} = D + D_1 + D_2.
$$
As to $I\!H^2(Y),$ the perversity condition is empty for $4-$chains 
since the singular locus is zero-dimensional. 
Hence $ \Delta_1, \Delta_2$ define intersection cohomology classes and in fact 
$$
I\!H^2 (Y) = \langle \Delta_1, \Delta_2 \rangle, 
\hbox { with the relation }
\Delta_{\infty} = \Delta_1 + \Delta_2.
$$
which is easily checked to hold on $Y-\{v \},$ hence on $Y.$

We have  $IC_X = \rat_X[3].$ The perverse cohomology complexes are
(cf. \ci{leiden}, \ci{htam}):
$^p\!H^{0} (f_{\bullet} IC_X ) = P^0_{\eta} = IC_Y;$
$^p\!H^{-1} (f_{\bullet} IC_X ) = P^{-1}_{\eta} = H_4(D)_v$
(skyscraper sheaf; in cohomological degree zero;
generated by the fundamental class of $D$);
$^p\!H^{1} (f_{\bullet} IC_X ) =  P^{-1}_{\eta} (-1) = H^4(D)_v.$

The Decomposition Theorem yields the existence of an isomorphism:
\begin{equation}
\label{isoqu}
\phi \; : \;  H_4(D)[1] \, \oplus  \, IC_Y \,\oplus \,
H^4(D)[-1] \; \simeq \;  f_{\bullet}IC_X.
\end{equation}

The resulting inclusion 
$\phi (I\!H^2 (Y) ) \subseteq H^2(X)$ depends on $\phi.$

Even with the choice $\phi = g_{\eta},$
the subspace 
$g_{\eta} (I\!H^2 (Y) )= f_{i,\eta} (I\!H^2 (Y) )
 \subseteq H^2(X)$ still depends on $\eta$ in a way we now make explicit.

The map $e: \, ^p\!H^{-1} (f_{\bullet} IC_X ) \simeq \,
^p\!H^{1} (f_{\bullet} IC_X ) $ 
is the map
\begin{equation}
\label{expe}
[D] \lorw e( [D]) \, = \, \{ [D] \lorw \eta \cdot D \cdot D \}, 
\end{equation}
where the product is in $H(X).$

By Corollary  \ref{ocond} we have:
$$
H^2(X) \; \supseteq  \; f_{i,\eta} ( I\!H^2(Y)     ) \; = \; 
\{ a \in H^2(X) \, | \; \eta \cdot D \cdot a =0 \}
$$
so that the dependence on $\eta$ is now transparent.

For example, set 
 $\eta:= mD_1+ D_2,$ 
 $m \in \rat^+.$ 
Then 
$$
g_{\eta} ( I\!H^2(Y) ) \;  = \; \langle D_{\infty}, -mD_1+D_2 \rangle  \; \subseteq 
H^2(X),
$$
$$
g_{\eta} ( \Delta_1 ) = D_1 + \frac{1}{m+1} D, \quad 
g_{\eta} ( \Delta_2 ) = D_2 + \frac{m}{m+1} D.
$$
The conclusion is that different choices of $\eta$ produce
different embeddings on $I\!H(Y)$ into $H(X).$

It is amusing to note the following.
For $m=0,$ $\eta = D_2$ is no longer $f-$ample, but
the relative Hard Lefschetz still holds since $D_2\cdot D \cdot D = -1
\neq 0$
and we have  $g_{D_2} (\Delta_1) = D_1+ D$ and $g_{D_2}(\Delta_2) = D_2.$
This decomposition can be seen as the one that arises canonically
by factoring  (in precisely one of the two possible ways!) 
$f: X \to Y$ through the small resolution of the quadric cone.

Note also that there is no isomorphism $\phi$
yielding $\phi (\Delta_i) = D_i,$ $i=1,2.$ This is because
$\phi ( \Delta_1 + \Delta_2) = \phi (\Delta_{\infty} ) = D_\infty
 \neq D_1+ D_2.$

All the embeddings of $I\!H^2(Y)$ into $H^2(X)$ are, in this example
where everything is of pure type $(1,1)$,
compatible with the Hodge structures.
In general, this is not true. Our main result, Theorem \ref{tuttoht},
is that we can arrange for this to be true in complete generality.

We conclude this section by remarking that the 
mechanism 
 in the proof  of Proposition \ref{so} becomes transparent in this example,
where $i=0.$
In fact, we  start with
{\em any} lift $y: IC_Y \to \tu{-1} f_{\bullet} IC_X=
f_{\bullet} IC_X$ of the natural map
$x : IC_Y \to  \, ^p\!H^0 (f_{\bullet} IC_X) \to
\tu{0} f_{\bullet} IC_X.$ 
The Snake Lemma allows to correct uniquely $y$
by adding to it a map
$IC_Y \to \, ^p\!H^{-1} (f_{\bullet} IC_X)[1]$
so that the resulting map
$f_{0,\eta} = \tau= x_{-1}: IC_Y
\to f_{\bullet} IC_X$
has the property that the composition
$$
IC_Y \stackrel{f_{0,\eta}}\lorw
 f_{\bullet} IC_X \stackrel{\eta}\lorw f_{\bullet} IC_X [2] 
\lorw 
\tu{1} f_{\bullet} \,  IC_X[2] =\,  ^p\!H^{1} (f_{\bullet} IC_X)[1]=
H^4(D)[1]
$$
is the zero map.
In hypercohomology, i.e. in $H(X),$ this translates
into the condition 
$$
\eta \cdot D \cdot f_{0,\eta} ( \Delta_i ) =0, 
\qquad i =1,2.
$$

\section{Applications}
\label{tdiht}
We give few applications of Theorem \ref{tuttoht}.
\subsection{The intersection pairing on $I\!H(Y).$}
\label{ippf}
\begin{tm}
\label{c5} Let $Y$ be a compact algebraic variety of dimension $n.$
For every $l \in \zed$ the intersection pairing
$$
d_Y\; :  \; I\!H^{n-l}(Y) \lorw 
I\!H^{n+l}(Y)^{\vee} (-n) 
$$
is an isomorphism of weight $(n-l)$ PHS.
\end{tm}

{\em Proof.}
Let $f: X \to Y$ be a projective resolution of the singularities
of $Y$ and $\eta$ be an $f-$ample line bundle
 on $X.$
There is the diagram
$$
\xymatrix{
f_{\bullet} IC_X \ar[r]^{f_{\bullet} d_X}   
& (  f_{\bullet} IC_X  )^{\vee}   \ar[d]^{  g_{\eta}^{\vee} }\\
IC_Y \ar[u]^{g_{\eta}} \ar[r]^{ d_Y} & IC_Y^{\vee}.
}
$$

The composition $g_{\eta}^{\vee} \circ f_{\bullet}d_X \circ g_{\eta}
\in {\rm Hom}(IC_Y, IC_Y^{\vee}) \simeq \rat$
coincides with $d_Y$ on the smooth locus
of $f$ on $Y$ and hence on the whole $Y;$
see Remark \ref{ranami}.

\n
It follows that the duality isomorphism  $d_Y$
can be exhibited as the composition
of maps of PHS
by virtue of Theorem \ref{tuttoht}:
$$
I\!H^{n-l}(Y) \stackrel{g_{\eta}}\lorw H^{n-l}(X) \stackrel{d_X}\simeq
H^{n+l}(X)^{\vee} (-n) \stackrel{ proj \circ g_{\eta}^{\vee}} \lorw
I\!H^{n+l}(Y)^{\vee} (-n).
$$  \blacksquare

\subsection{The map $H(Y) \to I\!H (Y)$} 
\label{s1}
Given a compact algebraic variety  $Y$ of dimension $n$
there is the natural map $a_Y: \rat_Y[n] \to IC_Y$
and the induced map in hypercohomology
$a_Y: H(Y) \to I\!H(Y).$

\bigskip
We freely employ the language and basic results of the theory of MHS
in \ci{ho3}.
The MHS on $H^l(Y)$ has weights $\leq l,$  i.e. $W_lH^l(Y) = H^l(Y).$
 In fact, for every
resolution of singularities $f: X \to Y,$ $\ke{f^*} = W_{l-1}H^l(Y).$
The quotient $H^l(Y)/W_{l-1} H^l(Y)$ is a PHS of weight $l.$

\begin{tm}
\label{alce} ({\bf   The natural map $ H(Y) \to I\!H(Y)$})
Let $Y$ be compact. The natural map 
$$
a_Y \, : \, H^{l}(Y) \lorw \ih{l}{Y}
$$
is a map of  MHS (the r.h.s.  is a PHS) and
$$
\ke{a_Y} \,  = \, 
\ke{f^{*}}  \, = \, W_{l -1} \, H^{l }(Y).
$$
\end{tm}
{\em Proof.} Let $f: X \to Y$ be a  {\em projective} resolution of
 singularities. Let $g_{\eta}$ be the 
isomorphism associated with some ample line bundle $\eta$ on $X$.
As in the proof of Theorem \ref{c5}, the formula (\ref{cit})
yields the commutative diagram
$$
\xymatrix{
\rat_{Y}[n] \ar[rr]^{adj(f)} \ar[dr]^{a_Y} & & \fxn{f}{X}{n}
 \\
 & IC_{Y} \ar[ur]^{g_{\eta}} &
 }
 $$
 with $g_{\eta}$ admitting a splitting $g'.$  Since $g_{\eta}$ is injective,
 $\ke{a_Y} = \ke{f^*}$  and the second conclusion follows.

\n
 We have $a_Y = g' \circ adj(f).$
The map induced in hypercohomology $a_Y = g' \circ f^*$ is
the composition of the map $f^*$ of MHS (Deligne's
theory of MHS) and of the splitting $g'$
of PHS (Theorem \ref{tuttoht}).

\n
If $f$ is not projective, then the second assertion follows by considering
a Chow envelope
$f': X' \to X \to Y$ with $f'$ projective and then by recalling
 that Deligne's theory of MHS ensures that
$\ke{f^*}=  \ke{ {f'}^*}.$
\blacksquare

\subsection{Projectors and Hodge classes}
\label{s4}
Let $f :X \to Y$ be a projective morphism
of proper varieties and $\eta$ be an $f-$ample line bundle on $X.$

Let $I\!H^l_i(X) = H \oplus H'$
be a direct sum decomposition into SHS.

Using the decomposition
isomorphism $g_{\eta}$  and the associated projectors
we obtain the composition
$$
p_H \; : \; I\!H^l (X) \lorw H \lorw
I\!H^l (X)
$$
which is a projector, i.e. $p^2 =p$ in the algebra
$$
{\rm End } ( I\!H(X) ) \; =\; I\!H(X) \otimes
I\!H(X)^{\vee} \; \simeq \; 
I\!H(X) \otimes I\!H(X)= I\!H (X \times X),
$$
where the middle isomorphisms stems from the nondegenerate intersection pairing
of Theorem \ref{c5}.

\begin{tm}
\label{prcon}
$$
p_H \; \in \; I\!H^{n,n}_{\rat} (X \times X).
$$
\end{tm}
{\em Proof.}
The proof is identical to the analogous one
to be found in \ci{leiden} for the case when $X$ is smooth.
The only missing piece is Theorem \ref{c5}.
\blacksquare.

\medskip
It is natural to ask whether the classes $p_H$ of Theorem \ref{prcon}
are algebraic, i.e. representable in terms
of admissible geometric chains arising from algebraic subvarieties.

If $X$ is smooth, then  this amounts to ask whether
these classes are in
$$\im{ (A_n(X \times X) \to H^{2n}(X \times X))}.$$
This takes one to the realm of the Standard Conjectures for algebraic cycles
and we have nothing to say in this direction, except for very special, 
yet non-trivial  cases.
In \ci{decmigmot}, we have shown that for proper semismall maps from
a smooth space, for every $H,$ one can find  canonical {\em algebraic }  
projectors $c'$ 
of type $(n,n).$  The key point is that $\dim{X \times_Y X} =n.$
The paper \ci{chowhilb}, dealt with the case of Hilbert schemes
of points on surfaces.
In  \ci{leiden}, we have shown that the same can be done
for the resolution of singularities of a threefold. The key point there is the use of the
$(1,1)-$Theorem.

\smallskip
We remark that if $H = I\!H_i^l(X)$ is a summand
as in (\ref{byst}),  a summand of the $(\eta, L)-$dec\-om\-position
of \ci{htam}, or any intersection of the two, then
it can be shown that the cycles $p_H$ are {\em absolute Hodge classes}
in the sense of  \ci{ahc}. We plan to pursue this aspect and some of its consequences 
in a forthcoming paper.

\subsection{Induced morphisms in intersection cohomology}
\label{tmic}

Intersection cohomology is not functorial in the ``space'' variable.
However, the paper \ci{cinque} constructs, for every proper map
$f: X \to Y,$  non-canonical morphisms
$I\!H^l(Y) \to I\!H^l(X).$

If $f$ is surjective, these morphisms
stem from the Decomposition Theorem and are splitting injections.
We now show how to choose them so that
they are map of PHS.

\begin{tm}
\label{imi}
Let $f: X^n \to Y^m$ be a projective, surjective map of compact varieties
of the indicated dimensions, $\rho: =n-m$
and $\eta$ be an $f-$ample line bundle on $X.$
Then there are a canonical splitting injection
\begin{equation}
\label{gamma}
\gamma \; : \; IC_Y \lorw 
\, ^p\!H^{- \rho} ( f _{\bullet} IC_X ),
\end{equation}
and a commutative diagram of MHS
\begin{equation}
\label{cdomhs}
\xymatrix{
H^l(Y)   \ar[r]^{a_Y} \ar[d]^{f^*}   
&  I\!H^l(Y)  \ar[d]^{  g_{\eta} \circ \gamma }
\\
H^l(X)  \ar[r]^{ a_X} &  I\!H^l(X).
}
\end{equation}
\end{tm}
{\em Sketch of proof.}
Let $Y_m \subseteq Y$ be the dense stratum of a stratification for $f.$
The perverse sheaf $^p\!H^{- \rho} (f_{\bullet} IC_X)$ restricted to
$Y_m$ reduces to the shifted local system
${\cal L} [m],$ where ${\cal L}$ is the semisimple local system
of the $\rho-$dimensional irreducible components
of the typical fiber $f^{-1}(y).$

\n
Since $X$ is irreducible, the $\pi_1 (Y_m, y)-$invariants
${\cal L}^{\pi_1(Y_m)} \simeq \rat_{Y_m} \subseteq {\cal L},$
and the inclusion splits.

\n
The Decomposition Theorem implies that
$$
IC_Y \; = \; IC_Y (  {\cal L}^{\pi_1(Y_m)}   )
\; \subseteq \;  IC_Y (  {\cal L}   )
\; \subseteq \;
^p\!H^{- \rho} (f_{\bullet} IC_X)
$$
and that all the inclusions split canonically
This gives the map $\gamma$ and proves (\ref{gamma}).

\n
The diagram
$$
\xymatrix{
\rat_Y [m]   \ar[r]^{a_Y} \ar[d]^{adj(f)}   
&  IC_Y  \ar[d]^{  g_{\eta} \circ \gamma }
\\
f_{\bullet} \rat_X[m]  \ar[r]^{ a_X} & f_{\bullet} IC_X[- \rho]
}
$$
commutes in view of
the formula  (\ref{rana})
and of  Remark \ref{ranami} applied to $U= Y_m.$
The diagram (\ref{cdomhs}) is induced by it
by taking hypercohomology and is therefore commutative.

\n
The decomposition by strata  (\ref{byst}) and Theorem
\ref{tuttoht}
imply 
that
$g_{\eta} ( I\!H^{n+l}(Y, {\cal L}   ) )
 \; \subseteq \; I\!H^{n+l}(X)$
is a SHS.

\n
We are left with checking that
$I\!H(Y) \subseteq  I\!H(Y, {\cal L})$
is a SHS. Once this is done, we conclude
using Theorem \ref{alce}.

\n
Without loss of generality, we may assume that $X$ and $Y$ are normal.

\n
There is the Stein factorization $X \stackrel{f'}\to Y' \stackrel{\nu}\to
Y$ where $f'$ has connected fibers, $Y'$ is normal and $\nu$ is finite.
We have $IC_Y({\cal L}) = \nu_* IC_{Y'},$ so that
we may replace
$f: X \to Y, $ by $\nu: Y' \to Y, $ i.e. we may assume
that $f$ is finite.

\n
We have the commutative diagram
$$
\xymatrix{
Z \ar[rr]^{h} \ar[dr]^k & & X  \ar[dl]_f
 \\
 & Y
 }
 $$
 arising from the Galois closure of
$K(X)/K(Y).$ The maps $h$ and $k$ are finite Galois.

\n
The inclusions $I\!H(Y) \subseteq I\!H(X) \subseteq I\!H(Z)$
imply
that if we can prove
the wanted conclusion for a Galois map, then it
will follow for $f.$ This means that we may assume
that $f$ is Galois with finite Galois group $G.$

\n
The group $G$ acts on the PHS $I\!H(X)$ by automorphisms
of PHS: take a $G-$equivariant resolution of
the singularities $p: X' \to X,$ a $G-$invariant $p-$ample
line bundle $\eta'$ on $X'$ and use Theorem \ref{tuttoht}.

\n
It follows that the $g-$invariants
$I\!H(Y) = I\!H(X)^G \subseteq  I\!H(X)$ form  a SHS.
\blacksquare


Authors' addresses:

\smallskip
\n
Mark Andrea A. de Cataldo,
Department of Mathematics,
SUNY at Stony Brook,
Stony Brook,  NY 11794, USA. \quad 
e-mail: {\em mde@math.sunysb.edu}

\smallskip
\n
Luca Migliorini,
Dipartimento di Matematica, Universit\`a di Bologna,
Piazza di Porta S. Donato 5,
40126 Bologna,  ITALY. \quad
e-mail: {\em migliori@dm.unibo.it}


\begin{thebibliography}{99}

\bibitem{cinque}{G.Barthel J-P.Brasselet, K-H.Fieseler, O.Gabber, L.Kaup,
``Rel\`evement de cycles alg\'ebriques et homomorphismes associ\'es en 
homologie d'intersection,'' Ann. of Math. {\bf 141}(1995),  147-179.}
\bibitem{bbd}{A.A. Beilinson, J.N. Bernstein, P. Deligne,
{\em Faisceaux pervers}, Ast\'erisque {\bf 100}, Pa\-ri\-s, Soc. Math. 
Fr. 1982.}


\bibitem{borel} A. Borel et al., {\em Intersection Cohomology}, 
Progress in Mathematics Vol. {\bf 50}, Birkh\"auser, Boston Basel 
Stuttgart 1984.

\bibitem{chowhilb} { M.A. de Cataldo, L. Migliorini, ``The Chow groups 
and the motive of the Hilbert schemes of surfaces,"
Journal of algebra 251 (2002), 824-848.}



\bibitem{demigsemi}{ M. de Cataldo, L. Migliorini, ``The 
Hard Lefschetz Theorem and the Topology of semismall maps,"      
Ann. Sci. \'Ecole Norm. Sup. (4) {\bf 35} (2002), no. 5, 759--772.}

\bibitem{decmigmot}{ M. de Cataldo, L. Migliorini, ``The Chow Motive 
of semismall resolutions,''
Math.Res.Lett. 11, (2004), 151-170.}

\bibitem{htam}{M. de Cataldo, L. Migliorini, ``The 
Hodge Theory of Algebraic maps,''
Ann.Sci. \'Ecole Norm. Sup. {\bf 38} (2005)no.5, 693-750.}

\bibitem{leiden}{M. de Cataldo, L. Migliorini, ``
Intersection forms , topology of maps and motivic decomposition for resolutions of threefolds,'' to appear in the Proceedings of the conference in honour of J.Murre, Leiden 2004.}


\bibitem{dess} P. Deligne, ``Th\'eor\`eme de Lefschetz et crit\`eres 
de d\'eg\'en\'erescence
de suites spectrales,'' Publ.Math. IHES {\bf 35} (1969), 107-126.

\bibitem{shockwave} P. Deligne, ``D\'ecompositions dans la cat\'egorie 
D\'eriv\'ee'', 
Motives (Seattle, WA, 1991), 115--128, Proc. Sympos. Pure Math., {\bf 
55},Part 1, Amer. Math. Soc., Providence, RI, 1994.
 

\bibitem{ho3}  P. Deligne, ``Th\'eorie de Hodge, III,'' Publ.Math. 
IHES {\bf 44} (1974), 5-78.

\bibitem{ahc}{P.Deligne, ``Hodge cycles on abelian varieties,'' in
P.Deligne, J.S. Milne, A. Ogus, K.Shih, {\em Hodge cycles, motives and Shimura varieties }, LNM 900, Springer, 1982.}


\bibitem{go-ma2} M. Goresky, R. MacPherson, ``Intersection homology 
II,''
Inv. Math. {\bf 71} (1983), 77-129. 


\bibitem{mac83}{ R. MacPherson}{ ``Global Questions in the Topology
of Singular Spaces,'' Proc. of the I.C.M.,  213--235 (1983), Warszawa.}


 \bibitem{samhp}{    
 M. Saito, ``Modules de Hodge polarisables,''
 Publ. Res. Inst. Math. Sci. 
{\bf 24} (1988), no.~6, 849--995 (1989).}



\end{thebibliography}
\end{document}